\documentclass[12pt, a4paper, reqno, english]{amsart}
\usepackage{amssymb,amsmath,amsthm}
\numberwithin{equation}{section}
\usepackage[margin=1in]{geometry}
\usepackage[colorlinks=true, linkcolor=blue, citecolor=green, urlcolor=black]{hyperref}

\newcommand{\C}{\mathbb{C}}
\newcommand{\N}{\mathbb{N}}
\newcommand{\Z}{\mathbb{Z}}
\newcommand{\U}{\mathbb{U}}
\newcommand{\e}{\mathrm e}

\newtheorem{thm}{Theorem}[section]
\newtheorem{prop}{Proposition}[section]
\newtheorem{lem}{Lemma}[section]

\newtheorem{cor}{Corollary}[section]
\newtheorem{conjecture}{Conjecture}[section]

\theoremstyle{definition}

\usepackage{color}
\definecolor{red}{rgb}{1,0,0}
\definecolor{blue}{rgb}{.2,.6,.75}
\definecolor{green}{rgb}{.4,.7,.4}

\renewcommand{\leq}{\leqslant}
\renewcommand{\geq}{\geqslant}
\renewcommand{\d}{\mathrm{d}}

\begin{document}

\title[Logarithmic exponential sums]{Logarithmic exponential sums with multiplicative coefficients}

\author{Haozhe Gou}
\date{}

\address{School of Mathematics, Shandong University, Jinan 250100, China}
\address{D\'epartement de math\'ematiques et de statistique, Universit\'e de Montr\'eal, C.P.~6128, succ.~Centre-ville, Montr\'eal, QC H3C~3J7, Canada}
\email{haozhegou@gmail.com}


\begin{abstract}
Let $f$ be a $1$-bounded multiplicative function, and suppose that
$|\alpha-a/q|\leq 1/(qx)$, where $(a,q)=1$ and $q\leq x$. We prove that
\[
 \left|\sum_{n\leq x}\frac{f(n)}n\e(n\alpha)\right|
 \leq
 \left(\frac4{\sqrt6}+o(1)\right)\frac{\log x}{\sqrt q}
 +O(\log q).
\]
Here the constant $4/\sqrt{6}$ is optimal. 
This settles  two conjectures of Granville and Lamzouri recently. If $f$ is completely multiplicative, the constant improves to $1$.
The main new ingredient is a two-prime Ramar\'e decomposition identity.
\end{abstract}

\maketitle

\section{Introduction}

Let \(
 \U:=\{z\in\C:|z|\leq1\}\) and  $f:\N\to\U$ be multiplicative functions.
Exponential sums with multiplicative coefficients have been studied
extensively and have many applications in analytic number theory. A foundational  result was established by
Montgomery and Vaughan \cite{MontgomeryVaughan1977}, who obtained strong
bounds for such sums in terms of rational approximation to $\alpha$. 
The subject has since been developed further in several directions; see,
for example,
\cite{GranvilleSoundararajan2007,Goldmakher2012,
BoberGoldmakherGranvilleKoukoulopoulos2018,
LamzouriMangerel2022,BretecheGranville2022,
GranvilleLamzouri2026}
and the references therein.

A folklore conjecture predicts that, for each fixed $\varepsilon>0$, if
$|\alpha-a/q|\leq1/(qx)$, $(a,q)=1$, and
$q\leq x/(\log x)^{2+\varepsilon}$, then\footnote{We will prove this conjecture in forthcoming work.}
\begin{equation*}
 \sum_{n\leq x}f(n)\e(n\alpha)
 \ll \frac{x}{\log x}+\frac{x}{\sqrt q}.
\end{equation*}
In this paper, we are concerned with the logarithmically weighted analogue
\begin{equation}\label{eq:L-definition}
 L_f(x,\alpha):=\sum_{n\leq x}\frac{f(n)}n\e(n\alpha).
\end{equation}
Besides its relation to the unweighted sums through partial summation,
this logarithmic variant also arises naturally for Dirichlet characters in
P\'olya's Fourier expansion for character sums. 
Granville and Lamzouri recently proved, among other results, that
\begin{equation}\label{eq:GL-bound}
 |L_f(x,\alpha)|
 \ll \frac{\log x}{\sqrt{\phi(q)}}+\log q
\end{equation}
holds whenever $|\alpha-a/q|\le 1/(qx)$ with $(a,q)=1$ and $q\le x$.  Their examples show that the natural  size of the first term above should be $(\log x)/\sqrt q$, while an initial range $n\le q^2$ (say) can contribute to  the size $\log q$.  They formulated the
following two conjectures.

\begin{conjecture}[Granville--Lamzouri {\cite[(1.5)]{GranvilleLamzouri2026}}]
\label{conj:GL}
Let $f:\N\to\U$ be multiplicative and suppose that
$|\alpha-a/q|\leq1/(qx)$, where $(a,q)=1$ and $1\leq q\leq x$.  Then
\begin{equation}\label{eq:GL-conjecture}
 |L_f(x,\alpha)|\ll \frac{\log x}{\sqrt q}+\log q.
\end{equation}
\end{conjecture}

For the more precise form, 
we write $A(x)\lesssim B(x)$ if
$|A(x)|\leq(1+o_x(1))B(x)$ as $x\to\infty$.

\begin{conjecture}[Granville--Lamzouri {\cite[(1.6)]{GranvilleLamzouri2026}}]
\label{conj:GL-sharp}
Under the hypotheses of Conjecture~\ref{conj:GL}, if
$q\leq(\log x/\log\log x)^2$, then
\begin{equation}\label{eq:GL-sharp-conjecture}
 |L_f(x,\alpha)|\lesssim
 \frac4{\sqrt6}\frac{\log x}{\sqrt q}.
\end{equation}
\end{conjecture}

On the major arcs $q\le (\log x)^{2+\varepsilon}$, by adapting the
pretentious approach of Granville, Harper and
Soundararajan \cite{GranvilleHarperSoundararajan2019} to the logarithmically weighted setting, Granville and Lamzouri obtained
an asymptotic formula for $L_f(x,\alpha)$ in terms of some primitive
characters; see \cite[Theorem~1.6, Corollary~1.3]{GranvilleLamzouri2026}. As a consequence, they proved Conjecture~\ref{conj:GL-sharp} in the range
\(
 q\ll (\log x)^2\exp(-c\sqrt{\log\log x})
\)
for some constant $c>0$. 
Together with their general estimate \eqref{eq:GL-bound}, this leaves Conjecture~\ref{conj:GL} open only in the interval
\[
 (\log x)^2\exp(-c\sqrt{\log\log x})
 <q<
 \left(\frac{\log x}{\log\log x}\right)^2\log\log\log x.
\]
Therefore for the proof of our result, it will  be convenient to restrict to the range
\[
 \log x\leq q\leq \log^2 x.
\]

Our main result is the following.

\begin{thm}\label{thm:main}
Let $f:\N\to\U$ be multiplicative, and suppose that
$|\alpha-a/q|\leq 1/(qx)$, $(a,q)=1$, and $1\leq q\leq x$. Then
\begin{equation}\label{eq:main}
 |L_f(x,\alpha)|
 \leq
 \left(\frac4{\sqrt6}+o(1)\right)
 \frac{\log x}{\sqrt q}
 +O(\log q),
\end{equation}
If $f$ is completely multiplicative, then the constant $4/\sqrt{6}$ improves to $1$. Both constants are best possible in their respective classes.
\end{thm}

Theorem~\ref{thm:main} immediately implies Conjecture~\ref{conj:GL}. It
also gives the sharp form of Conjecture~\ref{conj:GL-sharp} whenever the
secondary term is negligible.

\begin{cor}\label{cor:GL-conjecture}
Conjecture~\ref{conj:GL} holds.
\end{cor}

\begin{cor}\label{cor:subcritical}
Let $f:\N\to\U$ be multiplicative, and suppose that
$|\alpha-a/q|\leq1/(qx)$, $(a,q)=1$, and
$q=o((\log x)^2/(\log\log x)^2)$. Then
\begin{equation}\label{eq:subcritical}
 |L_f(x,\alpha)|
 \leq
 \left(\frac4{\sqrt6}+o(1)\right)
 \frac{\log x}{\sqrt q}.
\end{equation}
If $f$ is completely multiplicative, then $4/\sqrt{6}$ improves to 1.
Both constants are best possible in their respective classes. 
\end{cor}

We can immediately show that Conjecture~\ref{conj:GL-sharp} cannot hold
 throughout the full range \(q\leq({\log x}/{(\log\log x)})^2
\)
with the same constant $4/\sqrt{6}$. 
Fix $0<c\leq1$ and take
\(
 q\sim
 c({\log x}/{(\log\log x)})^2,
  f\equiv1\) and \( \alpha=1/q.
\)
Then
\[
 L_f(x,\alpha)=\sum_{n\le x}\frac{e(n/q)}{n}
 =-\operatorname{Log}(1-\e(1/q))+O(q/x),
\]
where $\operatorname{Log}$ denotes the principal branch, and hence
\[
 |L_1(x,1/q)|
 =\log q+O(1)
 =\left(2\sqrt c+o(1)\right)\frac{\log x}{\sqrt q}.
\]
Thus the coefficient $4/\sqrt6$ fails  whenever
$c>2/3$, and the coefficient $1$   fails whenever $c>1/4$. 
Thus the range $q=o((\log x/\log\log x)^2)$ in Corollary~\ref{cor:subcritical} is
essentially best possible, in the sense that it cannot in general be
extended to the full  range
$q\asymp(\log x/\log\log x)^2$.


\subsection{Structure of the proof}

We first reduce both cases in Theorem~\ref{thm:main} to the range
$\log x\leq q\leq(\log x)^2$. For multiplicative functions, this follows
directly as noted above. Though Granville and Lamzouri also treat the
completely multiplicative case
\cite[Theorem~1.6 and Corollary~7.1]{GranvilleLamzouri2026}, to obtain the
sharp coefficient $1$ needed here, we need to combine their major arc asymptotic
formula with the corresponding local factors of the Euler product. We record this short
calculation in Section~\ref{sec:reductions}.

We now describe the main argument. Throughout this paper, let $y:=q^2$. Integers with at most
one prime factor exceeding $y$ give a lower order contribution; see Section \ref{sec:sparse}. For the
remaining integers we use a two-prime Ramar\'e identity (see Section \ref{sec:ramare-gauss}) to expose two
distinct large prime factors and write
\[
 n=dcpr,\qquad d\mid q^\infty,\qquad (cpr,q)=1,
\]
where $d\mid q^\infty$ means that all prime factors of $d$ divide $q$. 
Multiplicativity then separates the four factors. Fourier inversion on
$(\Z/q\Z)^\times$ transforms the additive phase into a character sum
involving the Gauss sum $G_q(\chi;d)$. The two marked prime variables are then estimated by a second moment over
Dirichlet characters, using the distribution of primes in reduced residue classes. The contribution of the factor $d$, which is supported on the
primes dividing $q$, is handled separately. 

The $d$-sum involving the Gauss sums can be estimated prime by prime by the Chinese remainder theorem. For multiplicative functions $f$ we obtain
\[
 \Bigg|
 \sum_{d\mid q^\infty}\frac{f(d)G_q(\chi;d)}{d}
 \Bigg|
 \leq
 \sum_{d\mid q^\infty}\frac{|G_q(\chi;d)|}{d}
 \leq \frac4{\sqrt6}\sqrt q.
\]
At a prime power $p^e\Vert q$, the local factor is at most $p^{e/2}$  except when the local character $\chi \pmod{p^e}$ is principal and $p^e=2$ or $3$.
These two exceptional factors contribute $2/\sqrt2$ and $2/\sqrt3$,
respectively. This is where the constant $4/\sqrt6$ comes from.
For completely multiplicative $f$, using $f(p^v)=f(p)^v$, we get
\[
 \Bigg|
 \sum_{d\mid q^\infty}\frac{f(d)G_q(\chi;d)}d
 \Bigg|
 \leq\sqrt q.
\]
 Apart from this distinction, the rest of the argument is the same in the two cases.
 
Finally, 
the weights from the two-prime Ramar\'e identity recombine into the harmonic sum over integers coprime
to $q$, with no further loss in the constants $4/\sqrt6$ and $1$. We estimate the
contribution from integers with at least two large prime factors in Sections~\ref{sec:two-large}, and then complete the proof of Theorem~\ref{thm:main} in Sections~\ref{sec:completion}.


\subsection*{Acknowledgement}
The author thanks Professors Andrew Granville, Youness Lamzouri and Jianya Liu for their helpful discussions and valuable feedback.

\section{Preliminary reductions}\label{sec:reductions}

We first replace $\alpha$ by $a/q$. Since
$|\e(u)-\e(v)|\leq 2\pi|u-v|$, we have
\begin{equation}\label{eq:phase-error}
 |L_f(x,\alpha)-L_f(x,a/q)|\leq \frac{2\pi}{q}.
\end{equation}
Thus it is enough to work at the rational point $a/q$. 

For general multiplicative functions, the range $q<\log x$ is already
covered by \cite[Corollary~1.3]{GranvilleLamzouri2026}. 
For completely multiplicative functions we treat as follows.

\begin{lem}\label{lem:complete-small-q}
Let $f:\N\to\U$ be completely multiplicative. Let $2\leq q<\log x$ and $|\alpha-a/q|\leq1/(qx)$ with $(a,q)=1$, we have
\[
 |L_f(x,\alpha)|
 \leq (1+o(1))\frac{\log x}{\sqrt q}.
\]
\end{lem}

\begin{proof}
By \cite[Corollary~7.1]{GranvilleLamzouri2026}, there is a primitive
character $\chi\pmod\ell$, with $\ell\mid q$, such that
\begin{equation}\label{eq:GL-complete}
 L_f(x,\alpha)
 =
 \frac{\overline{\chi}(a)\kappa(q/\ell)\tau(\chi)}{\phi(q)}
 \sum_{n\leq x}\frac{f(n)\overline{\chi}(n)}n
 +O\bigg(\frac{(\log x)^{2/3+o(1)}}{\sqrt q}+e^{O(\sqrt{\log\log x})}\bigg),
\end{equation}
where $f=\kappa*\chi$. The error term is clearly 
$o((\log x)/\sqrt q)$ since $q<\log x$.

For every $n$ with $\chi(n)\neq0$, writing $n=dm$, where all prime factors of $d$ divide $q$ but not $\ell$,
and $(m,q)=1$, we obtain
\[
\begin{aligned}
 \sum_{n\leq x}\frac{f(n)\overline{\chi(n)}}n
 &=
 \sum_{\substack{d\leq x\\
                  p\mid d\Rightarrow p\mid q,\ p\nmid\ell}}
 \frac{f(d)\overline{\chi(d)}}d
 \sum_{\substack{m\leq x/d\\(m,q)=1}}
 \frac{f(m)\overline{\chi(m)}}m\\
 & 
 = \prod_{\substack{p\mid q\\p\nmid\ell}}
 \bigg(1-\frac{f(p)\overline{\chi(p)}}p\bigg)^{-1}
 \sum_{\substack{m\leq x\\(m,q)=1}}
 \frac{f(m)\overline{\chi(m)}}m 
 +O\left(\log^2q\right),
\end{aligned}
\]
where we replace $m\le x/d$ by $m\le x$ in the inner sum of the first line, at a cost
\(
 \sum_{x/d<m\leq x}m^{-1}
 \ll \log d\), and then extend the $d$-sum to infinity. A crude estimate gives  the error
 \[
 \sum_{\substack{d\ge1\\p\mid d\Rightarrow p\mid q}}
 \frac{\log d}{d}
 =
 \prod_{p\mid q}\bigg(1-\frac1p\bigg)^{-1}
 \sum_{p\mid q}\frac{\log p}{p-1}
 \ll
 \frac q{\phi(q)}\log q
 \ll \log^2 q,
 \]
while the contribution from $d>x$ is o(1) by Rankin's trick.

Moreover, by M\"obius inversion,
\[
 \Bigg|
 \sum_{\substack{m\leq x\\(m,q)=1}}
 \frac{f(m)\overline{\chi(m)}}m
 \Bigg|
 \leq
 \sum_{\substack{m\leq x\\(m,q)=1}}\frac1m
 =\sum_{r\le x,r|q}\frac{\mu(r)}{r}\sum_{n\le x/r}\frac{1}{n}
 =
 \frac{\phi(q)}q\log x+O(\log^2 q).
\]
Thus, using $|\tau(\chi)|=\sqrt\ell$, it remains to show that
\begin{equation}\label{eq:complete-local-bound}
 \frac{\sqrt\ell}{q}|\kappa(q/\ell)|
 \prod_{\substack{p\mid q\\p\nmid\ell}}
 \Bigg|1-\frac{f(p)\overline{\chi(p)}}p\Bigg|^{-1}
 \leq \frac1{\sqrt q}.
\end{equation}

We verify this prime by prime. Write $p^e\Vert q$ and $p^j\Vert\ell$,
where $0\leq j\leq e$. If $j\geq1$, then $\chi(p)=0$ and
$\kappa(p^{e-j})=f(p)^{e-j}$. Hence the corresponding local factor on
the left-hand side of \eqref{eq:complete-local-bound} is at most
$p^{j/2-e}\leq p^{-e/2}$.
If $j=0$, put $z:=f(p)\overline{\chi(p)}$. Since $f=\kappa*\chi$ and
$f$ is completely multiplicative, we have
\[
 \kappa(p^e)=f(p)^e-\chi(p)f(p)^{e-1}\implies |\kappa(p^e)|=|f(p)^{e-1}||1-z|\le|1-z|.
\]
The corresponding local factor on
the left-hand side of \eqref{eq:complete-local-bound} is therefore at most
\[
 \frac1{p^e}\frac{|1-z|}{|1-z/p|}\le p^{-e+1/2}\le p^{-e/2},
\]
where we used $|1-z|\leq\sqrt p |1-z/p|$, since $|z|\leq1$ and 
\(
 p|1-z/p|^2-|1-z|^2
 =(p-1)\left(1-{|z|^2}/{p}\right)\geq0.
\)
Multiplying over the local factor at $p^e\Vert q$ proves
\eqref{eq:complete-local-bound}.

Substituting \eqref{eq:complete-local-bound}  into \eqref{eq:GL-complete} we get
\[
 |L_f(x,\alpha)|
 \leq (1+o(1))\frac{\log x}{\sqrt q}.
\]
The two $O(\log^2 q)$ errors above contribute, after substitution into
\eqref{eq:GL-complete}, at most
\[
\begin{aligned}
 \log^2 q\left(
 \frac{\sqrt{\ell}}{\phi(q)}|\kappa(q/\ell)|
 +\frac{\sqrt q}{\phi(q)}
 \right)
 &\ll
 \frac{\sqrt q}{\phi(q)}
 \bigl(2^{\omega(q)}+1\bigr)\log^2 q=
 o\left(\frac{\log x}{\sqrt q}\right),
\end{aligned}
\]
where the second term follows from \eqref{eq:complete-local-bound} and the first term comes from
\(
 |\kappa(a)|
\le
 2^{\omega(a)}.
\)
This proves the lemma.
\end{proof}

At the other end, if $q>(\log x)^2$, then \eqref{eq:GL-bound} and the
standard estimate $q/\phi(q)\ll\log\log(3q)$ give us
\(
 {\log x}/{\sqrt{\phi(q)}}\ll\log q.
\)
Thus both assertions of Theorem~\ref{thm:main} follow in this range.

For the rest of the proof we may therefore assume
\( \log x\leq q\leq(\log x)^2.
\)
In particular,
\begin{equation}\label{eq:direct-facts}
 \log q\asymp\log\log x,
 \qquad
 \frac q{\phi(q)}\ll\log\log\log x.
\end{equation}
The integers now divisible by $p^2$ for some prime $p>y$
have total harmonic weight
\[
 \ll(\log x)\sum_{p>y}\frac1{p^2}
 \ll\frac{\log x}{y}
 =o\left(\frac{\log x}{\sqrt q}\right).
\]
After deleting them, let $L_0$, $L_1$, and $L_{\geq2}$ denote the
parts of $L_f(x,a/q)$ for which $\omega_y(n)$ is respectively $0$, $1$,
and at least $2$. Thus we have
\begin{equation}\label{eq:decomposition}
 L_f(x,a/q)
 =
 L_0+L_1+L_{\geq2}
 +o\left(\frac{\log x}{\sqrt q}\right).
\end{equation}

\section{Integers with at most one large prime factor}\label{sec:sparse}

We show that the contribution from integers with at most one prime
factor exceeding $y=q^2$ is negligible. 

\begin{prop}\label{prop:sparse}
In the range $\log x\leq q\leq(\log x)^2$, we have
\[
 |L_0|+|L_1|
 \ll
 \log\log x
 +o\bigg(\frac{\log x}{\sqrt q}\bigg).
\]
\end{prop}

\begin{proof}
The contribution of $L_0$ is immediate from Mertens' theorem:
\[
 |L_0|=\bigg|
 \sum_{\substack{n\leq x\\P^+(n)\leq y}}
 \frac{f(n)}n\e(an/q)
 \bigg|
 \leq
 \prod_{p\leq y}\bigg(1-\frac1p\bigg)^{-1}
 \ll\log y
 \ll\log\log x.
\]

We now turn to $L_1$. After removing the terms divisible by $p^2$ for
some $p>y$, as in the previous section, every remaining integer counted
by $L_1$ can be written uniquely as $n=pm$, where $p>y$ and
$P^+(m)\leq y$. Hence
\begin{equation}\label{eq:L1-formula}
 L_1
 =
 \sum_{p>y}\frac{f(p)}p
 \sum_{\substack{m\leq x/p\\P^+(m)\leq y}}
 \frac{f(m)}m\e(amp/q).
\end{equation}

We first consider $m\leq\sqrt q$. We may replace $p\leq x/m$ by  $p\leq x/\sqrt q$, at a cost
\[
\begin{aligned}
 \sum_{m\leq\sqrt q}\frac1m
 \sum_{x/\sqrt q<p\leq x/m}\frac1p
 &\ll
 \log q\bigg(-\log\big(1-\frac{\log q}{2\log x}\big) +O\big(\frac1{\log x}\big) \bigg) \ll
 \frac{\log^2 q}{\log x}
 =o(1).
\end{aligned}
\]
Since $\sqrt q<y$, the condition $P^+(m)\leq y$ is automatic. 
Thus, after grouping the prime variable by its residue class modulo $q$,
the contribution from $m\leq\sqrt q$ is
\[
 \sum_{u\in(\Z/q\Z)^\times} A(u)V(u)
 \ll \|A\|_2\|V\|_2,
\]
by Cauchy--Schwarz, where
\[
 A(u):=
 \sum_{\substack{y<p\leq x/\sqrt q\\p\equiv u\pmod q}}
 \frac{f(p)}p,
 \qquad
 V(u):=
 \sum_{m\leq\sqrt q}\frac{f(m)}m\e(amu/q).
\]

For $\|A\|_2$, let $I=(P,2P]$ be a dyadic interval.  We write $$A(u)=\sum_I A_I(u),\quad A_I(u):=
\sum_{\substack{p\in I\\ p\equiv u\pmod q}}\frac{f(p)}p.$$ 
Since
$P\geq y=q^2$, $\log (P/q)\asymp\log P$, then by the Brun--Titchmarsh inequality and partial summation formula, we get
\[
 \max_{u\in(\Z/q\Z)^\times}|A_I(u)|\le \max_{u\in(\Z/q\Z)^\times}
 \sum_{\substack{p\in I\\p\equiv u\pmod q}}\frac1p
 \ll
 \frac1{\phi(q)}
 \sum_{p\in I}\frac1p.
\]
Hence by Minkowski's inequality, we have
\[
\begin{aligned}
 \|A\|_2
 &\leq \sum_I\|A_I\|_2
 \leq
 \sum_I
 \bigg(
 \max_u|A_I(u)|
 \sum_{u\in(\Z/q\Z)^\times}|A_I(u)|
 \bigg)^{1/2}\\
 &\leq
 \sum_I
 \bigg(
 \frac1{\phi(q)}
 \sum_{p\in I}\frac1p
 \sum_{p\in I}\frac1p
 \bigg)^{1/2}
 \ll
 \frac1{\sqrt{\phi(q)}}
 \sum_I\sum_{p\in I}\frac1p
 \ll
 \frac{\log\log x}{\sqrt{\phi(q)}}.
\end{aligned}
\]

For $\|V\|_2$, let $c_q(h):= \sum_{u\in(\Z/q\Z)^\times}\e(hu/q)$ denote the Ramanujan sum. For $m\neq n$, the standard bound $|c_q(h)|\le(q,h)$ gives
$|c_q(m-n)|\le(q,m-n)\le |m-n|<\sqrt q. $ Therefore
\[
\begin{aligned}
 \|V\|_2^2
 &=
 \sum_{u\in(\Z/q\Z)^\times}|V(u)|^2
 =
 \sum_{m,n\leq\sqrt q}
 \frac{f(m)\overline{f(n)}}{mn}\,c_q(m-n)\\
 &\ll
 \phi(q)\sum_{m\leq\sqrt q}\frac1{m^2}
 +
 \sqrt q\sum_{\substack{m,n\leq\sqrt q\\m\neq n}}\frac1{mn}\ll
 \phi(q)+\sqrt q \log^2 q
 \ll \phi(q).
\end{aligned}
\]
Combining this with the bound for $\|A\|_2$, the contribution from
$m\leq\sqrt q$ is \( \|A\|_2\|V\|_2\ll\log\log x.
\)

We now consider $m>\sqrt q$. Split the prime variable into dyadic
intervals $I=(P,2P]$. For $p\in I$, we replace $m\leq x/p$ by  $\sqrt q<m\leq x/(2P)$. The omitted boundary terms are
bounded by
\[
\begin{aligned}
 \sum_I\sum_{p\in I}\frac1p
 \sum_{x/(2P)<m\leq x/p}\frac1m
 &\ll
 \sum_I\sum_{p\in I}\frac1p
 \ll \log\log x.
\end{aligned}
\]

For a fixed $I$, grouping both variables by their residue classes
modulo $q$, the contribution from the interval $I$ is
\[
\begin{aligned}
\bigg|
 \sum_{u\in(\Z/q\Z)^\times} &A_I(u)\sum_{v\bmod q}
C_I(v)\e(auv/q) \bigg|
\le 
\|A_I\|_2 \Bigg(\sum_{u\bmod q}
\bigg|
\sum_{v\bmod q}C_I(v)\e(auv/q)
\bigg|^2\Bigg)^{1/2}\\
 &\le
 \|A_I\|_2 \Bigg(\sum_{v,w\bmod q}
C_I(v)\overline{C_I(w)}
\sum_{u\bmod q}\e(au(v-w)/q)\Bigg)^{1/2}
\le
 \sqrt q\,\|A_I\|_2\|C_I\|_2
 \end{aligned}
\]
where 
\[
C_I(v):=
 \sum_{\substack{\sqrt q<m\leq x/(2P)\\
                  P^+(m)\leq y\\
                  m\equiv v\pmod q}}
 \frac{f(m)}m,
\]
and
$$ \begin{aligned} \|C_I\|_2^2 &\leq \max_{v\bmod q}|C_I(v)| \sum_{v\bmod q}|C_I(v)|
\ll
\frac{\log x}{q}  \sum_{\substack{m\geq1\\P^+(m)\leq y}}\frac1m\ll  \frac{\log x\log\log x}{q}. \end{aligned} $$

Now, summing over the dyadic intervals $I$ and using the same estimate for $\|A_I\|_2$ as above, we obtain the contribution from $m>\sqrt{q}$ is
\[
\begin{aligned}
 \sqrt q\sum_I\|A_I\|_2\|C_I\|_2
 &\ll
 \sqrt q\sum_I
 \frac1{\sqrt{\phi(q)}}
 \bigg(\sum_{p\in I}\frac1p\bigg)
 \bigg(\frac{\log x\log\log x}{q}\bigg)^{1/2}\\
 &\ll
 (\log\log x)^{3/2}
 \sqrt{\frac{\log x}{\phi(q)}}
 =
 o\bigg(\frac{\log x}{\sqrt q}\bigg),
\end{aligned}
\]
where we used
\(
 \sum_I\sum_{p\in I}\frac1p\ll\log\log x
\)
and \eqref{eq:direct-facts}. Together with the boundary terms and the
range $m\leq\sqrt q$, this proves the proposition.
\end{proof}

\section{The Ramar\'e decomposition and Gauss sums}
\label{sec:ramare-gauss}

We now turn to $L_{\geq2}$.  A standard use of Ramar\'e's identity is to
extract one prime factor from an integer and thereby obtain a bilinear
decomposition.
The usual form (see, for example, \cite[Section~17.3]{FriedlanderIwaniec2010} and
\cite{MatomakiRadziwill2016}) is
\[
 \mathbf 1_{\omega_{\mathcal P}(n)\geq1}
 =
 \sum_{\substack{p\mid n\\p\in\mathcal P}}
 \frac1{\omega_{\mathcal P}(n)}
 =
 \sum_{\substack{pm=n\\p\in\mathcal P}}
 \frac1{\omega_{\mathcal P}(m)+\mathbf 1_{p\nmid m}},
\]
where  \( \omega_{\mathcal P}(n):=\#\{p\in\mathcal P:p\mid n\}
\) for some set of primes $ \mathcal P$.
For our purpose it is more convenient to extract two distinct primes
$p,r>y$ at the same time. This leaves two prime variables for the
character second moment used below.

We therefore use a slight reformulation of Ramar\'e's identity. The
elementary identity
\[
 \frac1{k(k-1)}
 =
 \int_0^1(1-t)t^{k-2}\,dt
 \qquad (k\geq2)
\]
replaces the combinatorial weight by the multiplicative weight
$t^{\omega_{\mathcal P}(\,\cdot\,)}$. This gives the following
two-prime form.


\begin{lem}[Two-prime Ramar\'e identity]\label{lem:ramare}
Suppose that no prime $p>y$ satisfies $p^2\mid n$ and that
$\omega_y(n)\geq2$. Then
\begin{equation}\label{eq:ramare}
 1=
 \int_0^1(1-t)
 \sum_{\substack{p,r>y,\ p\neq r\\pr\mid n}}
 t^{\omega_y(n/(pr))}\,\d t.
\end{equation}
\end{lem}

\begin{proof}
This follows directly from the identity we stated above by writting $k:=\omega_y(n)$.
\end{proof}

For each marked pair $p,r$, write
\[
 d:=\prod_{\ell\mid q}\ell^{v_\ell(n)},
 \qquad
 c:=\frac{n}{dpr}.
\]
Thus $d$ is the part of $n$ supported on the primes dividing $q$.
Since $p,r>y>q$, we have $(d,cpr)=(c,pr)=(p,r)=1$ and $(c,q)=1$.
Multiplicativity therefore gives $f(n)=f(d)f(c)f(p)f(r)$. Applying
Lemma~\ref{lem:ramare}, we get
\begin{equation}\label{eq:ramare-decomposition}
\begin{aligned}
 L_{\geq2}
 ={}&\int_0^1(1-t)
 \sum_{\substack{c\geq1,\ (c,q)=1\\
                  v_\ell(c)\leq1\ (\ell>y)}}
 \frac{f(c)t^{\omega_y(c)}}c
 \sum_{d\mid q^\infty}\frac{f(d)}d 
 \sum_{\substack{p,r>y,\ p\neq r\\
                  p\nmid c,\ r\nmid c\\dcpr\leq x}}
 \frac{f(p)f(r)}{pr}\e(adcpr/q)\,\d t.
\end{aligned}
\end{equation}

We now expand the additive character in Dirichlet characters modulo $q$.
For $\chi\pmod q$ and $n\in\Z$, write
\[
 \tau(\chi,n):=
 \sum_{u\bmod q}\chi(u)\e(nu/q)
\]
where $\tau(\chi,1)=\tau(\chi)$ is the Gauss sum, and if $(n,q)=1$, then
\(
 \tau(\chi,n)=\overline{\chi(n)}\,\tau(\chi).
\)
If $(m,q)=1$, then  by the orthogonality on $(\Z/q\Z)^\times$,
\begin{equation}\label{eq:character-expansion}
 \e(adm/q)
 =
 \frac1{\phi(q)}
 \sum_{\chi\bmod q}
 \tau(\overline{\chi},ad)\chi(m).
\end{equation}

\begin{lem}\label{lem:gauss}
For every $\chi\pmod q$, we have
\begin{equation}\label{eq:Gauss-l1}
 \sum_{d\mid q^\infty}
 \frac{|\tau(\overline{\chi},ad)|}{d}
 \leq
 \frac4{\sqrt6}\sqrt q.
\end{equation}
If $f$ is completely multiplicative, then
\begin{equation}\label{eq:Gauss-complete}
 \bigg|
 \sum_{d\mid q^\infty}
 \frac{f(d)\tau(\overline{\chi},ad)}{d}
 \bigg|
 \leq
 \sqrt q.
\end{equation}
Moreover, for every $\varepsilon>0$ and $X\geq1$,
\begin{equation}\label{eq:Gauss-tail}
 \sum_{\substack{d>X\\d\mid q^\infty}}
 \frac{|\tau(\overline{\chi},ad)|}{d}
 \ll_\varepsilon
 \frac{q^{1/2+\varepsilon}}{\sqrt X}.
\end{equation}
\end{lem}

\begin{proof}
By the Chinese remainder theorem, it is enough to estimate the
corresponding sums prime by prime. Let $q=p^e$, and suppose that $\chi$
is induced by a primitive character $\chi^\ast\pmod{p^j}$, where
$0\leq j\leq e$. If $j\geq1$, then
\[
\begin{aligned}
 \tau(\overline{\chi},ap^v)
 &=
 \begin{cases}
  0, & v\neq e-j,\\[2mm]
  p^{e-j}\chi^\ast(a)\tau(\overline{\chi^\ast}),
  & v=e-j,
 \end{cases}\\
 \sum_{v\geq0}
 \frac{|\tau(\overline{\chi},ap^v)|}{p^v}
 &=
 \frac{p^{e-j}|\tau(\overline{\chi^\ast})|}{p^{e-j}}
 =
 p^{j/2}
 \leq p^{e/2}.
\end{aligned}
\]
Here we used the standard formula
$|\tau(\chi^\ast)|=p^{j/2}$ for primitive characters.

If $j=0$, then $\chi$ is principal modulo $p^e$, and
$\tau(\overline{\chi},ap^v)=c_{p^e}(ap^v)$. Hence
\[
\begin{aligned}
 \sum_{v\geq0}\frac{|\tau(\overline{\chi},ap^v)|}{p^v}
 &=
 \frac{p^{e-1}}{p^{e-1}}
 +\phi(p^e)\sum_{v\geq e}\frac1{p^v}
 =2.
\end{aligned}
\]
Returning to general $q$, we get
\[
\begin{aligned}
 \sum_{d\mid q^\infty}\frac{|\tau(\overline{\chi},ad)|}{d}
 &\leq
 \sqrt q
 \prod_{\substack{p^e\Vert q\\
                   \chi\pmod{p^e}\ {\rm principal}}}
 \max\bigg(1,\frac2{p^{e/2}}\bigg)
 \leq
 \sqrt q\,\frac2{\sqrt2}\frac2{\sqrt3}
 =
 \frac4{\sqrt6}\sqrt q.
\end{aligned}
\]
This proves \eqref{eq:Gauss-l1}.

Suppose now that $f$ is completely multiplicative. By the Chinese
remainder theorem, the sum in \eqref{eq:Gauss-complete} factors prime
by prime, up to unit factors of modulus $1$. For $p^e\Vert q$, the
nonprincipal case contributes at most $p^{j/2}\leq p^{e/2}$. In the
principal case, after absorbing the unit factor into $f(p)$, the local
sum contributes, for some $|z|\le1$,
\[
 \bigg|z^{e-1}\frac{z-1}{1-z/p}\bigg|
 \leq\frac{|1-z|}{|1-z/p|}
 \leq\sqrt p
 \leq p^{e/2},
\]
by the same inequality $|1-z|\leq\sqrt p |1-z/p|$ we used in Lemma~\ref{lem:complete-small-q}.
Multiplying over $p^e\Vert q$ proves \eqref{eq:Gauss-complete}.

Finally,  we use Rankin's trick to bound
\[
 \sum_{\substack{d>X\\d\mid q^\infty}}
 \frac{|\tau(\overline{\chi},ad)|}{d}
 \leq
 X^{-1/2}
 \sum_{d\mid q^\infty}
 \frac{|\tau(\overline{\chi},ad)|}{\sqrt d}.
\]
The same calculation at prime power $p^e$ gives
\[
 \sum_{v\geq0}
 \frac{|\tau(\overline{\chi},ap^v)|}{p^{v/2}}
 \leq
 p^{e/2}\bigg(1+\frac2{\sqrt p}\bigg),
\]
and therefore \eqref{eq:Gauss-tail} follows by
\[
 \sum_{d\mid q^\infty}
 \frac{|\tau(\overline{\chi},ad)|}{\sqrt d}
 \ll
 \sqrt q\prod_{p\mid q}\bigg(1+\frac2{\sqrt p}\bigg)
 \ll_\varepsilon \sqrt q\,q^\varepsilon.
\]
\end{proof}

\section{Integers with at least two large prime factors}
\label{sec:two-large}

We now estimate $L_{\geq2}$. If one uses the Brun--Titchmarsh inequality
throughout, it is easy to obtain
\(
 L_{\geq2}\ll {(\log x)}/{\sqrt q}.
\)
However, to get the sharp constants $4/\sqrt6$ and $1$, we have to  make extra effort. 
Remember that the integer $n$ in this section is understood to satisfy
\(
 \omega_y(n)\geq2\) and \( p^2\nmid n \quad \text{for every prime } p>y.\) We will sometime omit these conditions.

\begin{prop}\label{prop:two-large}
In the range $\log x\leq q\leq(\log x)^2$, we have
\[
 |L_{\geq2}|
 \leq
 \bigg(\frac4{\sqrt6}+o(1)\bigg)
 \frac{\log x}{\sqrt q}+O(1).
\]
If $f$ is completely multiplicative, then
\[
 |L_{\geq2}|
 \leq
 (1+o(1))\frac{\log x}{\sqrt q}+O(1).
\]
\end{prop}

\begin{proof}
Split the prime variables into dyadic intervals
$I=(P,2P]$ and $J=(R,2R]$, we define
\[
  X:=\frac{x}{4cPR},\quad
P_{\chi,c}(I):=
 \sum_{\substack{p\in I\\p\nmid c}}
 \frac{f(p)\chi(p)}p \quad \text{and}
 \quad
 D_\chi(X):=
 \sum_{\substack{d\leq X\\d\mid q^\infty}}
 \frac{f(d)\tau(\overline{\chi},ad)}d.
\]
The choice of $X$ separates the condition $dcpr\leq x$ on each dyadic
block. 

If $d>X$, then
\(
 {x}/{4}<n=dcpr\leq x.
\)
After summing over $c,I,J$, taking absolute values and applying Lemma~\ref{lem:ramare}, the total
contribution of all such terms to $L_{\ge2}$ is at most
\[
\begin{aligned}
 &\int_0^1(1-t)
 \sum_{\substack{x/4<n\leq x}}
 \frac1n
 \sum_{\substack{p,r>y,\ p\neq r\\pr\mid n}}
 t^{\omega_y(n/(pr))}\,\d t
 \leq
 \sum_{x/4<n\leq x}\frac1n
 \ll1.
\end{aligned}
\]

For $d\leq X$, we have $dcpr\leq x$ for every
$p\in I$ and $r\in J$. After temporarily including the diagonal \(p=r\) when \(I=J\), whose contribution is acceptable, we have, by plugging
\eqref{eq:character-expansion} into \ref{eq:ramare-decomposition},

\begin{align}
 L_{\geq2}
 =&
 \int_0^1(1-t)
 \sum_{\substack{c\geq1,\ (c,q)=1\\
                  v_\ell(c)\leq1\;(\ell>y)}}
 \frac{f(c)t^{\omega_y(c)}}{c\phi(q)}
 \sum_{I,J}\sum_{\chi\bmod q}
 \chi(c)D_\chi(X)
 P_{\chi,c}(I)P_{\chi,c}(J)\d t
 +O(1)+o\bigg(\frac{\log x}{\sqrt q}\bigg)\notag\\
 ={}&
 \int_0^1 (1-t) \sum_{\substack{c\geq1,\ (c,q)=1\\
                  v_\ell(c)\leq1\;(\ell>y)}}t^{\omega_y(c)}
 \sum_{I,J}
 L_{\geq2,\mathrm{block}}(c;I,J)\,\d t
 +O(1) +o\bigg(\frac{\log x}{\sqrt q}\bigg),\label{eq:block}
\end{align}
where
\[
 L_{\geq2,\mathrm{block}}(c;I,J)
 :=
 \frac{f(c)}{c\phi(q)}
 \sum_{\chi\bmod q}
 \chi(c)D_\chi(X)
 P_{\chi,c}(I)P_{\chi,c}(J).
\]
By character orthogonality and the same Brun--Titchmarsh estimate used
in Section~\ref{sec:sparse},
\begin{align}\label{eq:prime-sec}
 \sum_{\chi\bmod q}|P_{\chi,c}(I)|^2
 &=
 \phi(q)
 \sum_{u\in(\Z/q\Z)^\times}
 \bigg|
 \sum_{\substack{p\in I,\ p\nmid c\\p\equiv u\pmod q}}
 \frac{f(p)}p
 \bigg|^2
 \ll
 \bigg(\sum_{p\in I}\frac1p\bigg)^2.
\end{align}
Then by \eqref{eq:Gauss-l1}, \eqref{eq:prime-sec}, and Cauchy--Schwarz, we have
\begin{equation}\label{eq:BT-block}
 |L_{\geq2,\mathrm{block}}(c;I,J)|
 \le
\frac{1}{c\phi(q)}
\max_{\chi}|D_\chi(X)|
\sum_{\chi\bmod q}
|P_{\chi,c}(I)P_{\chi,c}(J)|
 \ll
 \frac{\sqrt q}{c\phi(q)}
 \bigg(\sum_{p\in I}\frac1p\bigg)
 \bigg(\sum_{r\in J}\frac1r\bigg).
\end{equation}

We now fix $0<\delta<1$ and put \( z:=\exp((\log x)^\delta).\)
We first consider the part of $L_{\geq2}$ for which at least one of the two marked
primes lies between $y$ and $z$.
 By \eqref{eq:BT-block} and \eqref{eq:block}, and extending the sums, their total contribution is
\[
\begin{aligned}
 &\ll
 \frac{\sqrt q}{\phi(q)}
 \bigg(\sum_{y<p\leq z}\frac1p\bigg)
 \bigg(\sum_{y<p\leq x}\frac1p\bigg)
 \int_0^1(1-t)\sum_{\substack{c\leq x,\ (c,q)=1\\
                  v_\ell(c)\leq1\;(\ell>y)}}
 \frac{t^{\omega_y(c)}}c\,\d t
 +O(1)\\
 &\ll
 \frac{\sqrt q}{\phi(q)}
 \bigg(\sum_{y<p\leq z}\frac1p\bigg)
 \bigg(\sum_{y<p\leq x}\frac1p\bigg)
 \frac{\phi(q)}q\log x
 \int_0^1(1-t)
 \exp\bigg(-(1-t)\sum_{y<p\leq x}\frac1p\bigg)\,\d t
 +O(1)\\
 &\ll
 \frac{\log x}{\sqrt q}
 \bigg(\sum_{y<p\leq z}\frac{1}{p}\bigg)
      \bigg(\sum_{y<p\leq x}\frac{1}{p}\bigg)^{-1}
 +O(1)\\
 &\ll
 \frac{\delta+o(1)}{1+o(1)}
 \frac{\log x}{\sqrt q}
 +O(1).
\end{aligned}
\]
Here
\[
\begin{aligned}
 &\sum_{\substack{c\leq x,\ (c,q)=1\\
                  v_\ell(c)\leq1\;(\ell>y)}}
 \frac{t^{\omega_y(c)}}c\leq
 \prod_{\substack{p\leq y\\p\nmid q}}
 \bigg(1-\frac1p\bigg)^{-1}
 \prod_{y<p\leq x}\bigg(1+\frac tp\bigg)\\
 & =
 \prod_{\substack{p\leq x\\p\nmid q}}
 \bigg(1-\frac1p\bigg)^{-1}
 \prod_{y<p\leq x}
 \bigg(1-\frac1p\bigg)\bigg(1+\frac tp\bigg)
 \ll
 \frac{\phi(q)}q\log x\,
 \exp\bigg(-(1-t)\sum_{y<p\leq x}\frac1p\bigg).
\end{aligned}
\]

It remains to consider the part of $L_{\geq2}$ for which both marked
primes exceed $z$. Since $p\in(P,2P]$ and
$r\in(R,2R]$, we have $P,R\geq z/2$. Since \( q\leq(\log x)^2\leq(\log P)^{2/\delta},
\) we can apply the Siegel--Walfisz theorem to get
\begin{equation}\label{eq:prime-second-sharp}
\begin{aligned}
 \sum_{\chi\bmod q}|P_{\chi,c}(I)|^2
 &=
 \phi(q)
 \sum_{u\in(\Z/q\Z)^\times}
 \bigg|
 \sum_{\substack{p\in I,\ p\nmid c\\p\equiv u\pmod q}}
 \frac{f(p)}p
 \bigg|^2\\
 &\leq
 \phi(q)
 \max_{u\in(\Z/q\Z)^\times}
 \sum_{\substack{p\in I\\p\equiv u\pmod q}}\frac1p
 \sum_{p\in I}\frac1p
 \leq
 \bigl(1+o(1)\bigr)
 \bigg(\sum_{p\in I}\frac1p\bigg)^2.
\end{aligned}
\end{equation}
Suppose $X=x/(4cPR)<\log x$, then \( \frac{x}{4\log x}<cpr\leq x.\)
Using \eqref{eq:BT-block} and Lemma~\ref{lem:ramare}, the total
contribution of these blocks is
\[
\begin{aligned}
 &\ll
 \frac{\sqrt q}{\phi(q)}
 \sum_{x/(4\log x)<m\leq x}\frac1m\ll
 \frac{\sqrt q}{\phi(q)}\log\log x
 =
 o\bigg(\frac{\log x}{\sqrt q}\bigg),
\end{aligned}
\]
by \eqref{eq:direct-facts}, which is negligible.

We are left with $X\geq\log x$. For a general multiplicative function,
\eqref{eq:Gauss-l1} gives
\[
 |D_\chi(X)|
 \leq
 \frac4{\sqrt6}\sqrt q.
\]
If $f$ is completely multiplicative, then
\eqref{eq:Gauss-complete} and \eqref{eq:Gauss-tail} give, for $\varepsilon=1/8$ (say),
\[
\begin{aligned}
 |D_\chi(X)|
 &\leq
 \sqrt q+ O_\varepsilon\bigg(\frac{q^{1/2+\varepsilon}}{\sqrt X}\bigg)
 =
 (1+o_x(1))\sqrt q.
\end{aligned}
\]

Hence, by \eqref{eq:prime-second-sharp} and Cauchy--Schwarz,
\[
\begin{aligned}
 |L_{\geq2,\mathrm{block}}(c;I,J)|
 &\leq
 \frac{\max_\chi|D_\chi(X)|}{c\phi(q)}
 \bigg(\sum_{\chi\bmod q}|P_{\chi,c}(I)|^2\bigg)^{1/2}
 \bigg(\sum_{\chi\bmod q}|P_{\chi,c}(J)|^2\bigg)^{1/2}\\
 &\leq
 \bigg(\frac4{\sqrt6}+o_x(1)\bigg)
 \frac{\sqrt q}{c\phi(q)}
 \bigg(\sum_{p\in I}\frac1p\bigg)
 \bigg(\sum_{r\in J}\frac1r\bigg)
\end{aligned}
\]
in the multiplicative case. If $f$ is completely multiplicative, the same
estimate holds with $4/\sqrt6$ replaced by $1$.

To recombine the prime sums, we may drop the conditions $p\nmid c$ and
$r\nmid c$. 
The corresponding loss is
\[
 \ll
 (\log x)(1+\log\log x)
 \sum_{p>y}\frac1{p^2}
 =
 o\bigg((\log x)\frac{\phi(q)}q\bigg),
\]
which is again negligible  to $L_{\geq2}$.
For the remaining terms, put $m=cpr$ and write
\[
 K:=\omega_y(m),
 \qquad
 K_z:=\#\{p>z:p\mid m\}.
\]
For each fixed $m$, the marked primes can be chosen in
$K_z(K_z-1)$ ordered ways, while the two-prime Ramar\'e weight gives
\[
\begin{aligned}
 K_z(K_z-1)
 \int_0^1(1-t)t^{K-2}\,\d t
 &=
 \frac{K_z(K_z-1)}{K(K-1)}
 \leq1.
\end{aligned}
\]
Thus, in the multiplicative case, the contribution from $p,r>z$ and
$X\geq\log x$ is at most
\[
\begin{aligned}
 &\bigg(\frac4{\sqrt6}+o_x(1)\bigg)
 \frac{\sqrt q}{\phi(q)}
 \sum_{\substack{m\leq x\\(m,q)=1}}\frac1m
 =
 \bigg(\frac4{\sqrt6}+o_x(1)\bigg)
 \frac{\log x}{\sqrt q}
 +o\bigg(\frac{\log x}{\sqrt q}\bigg).
\end{aligned}
\]
For completely multiplicative $f$, the same argument gives
\(
 (1+o_x(1))\frac{\log x}{\sqrt q}.
\)

Combining this with the contribution from the range in which at least
one marked prime lies in $(y,z]$, we obtain
\[
 |L_{\geq2}|
 \leq
 \bigg(
 \frac4{\sqrt6}+O(\delta)+o_x(1)
 \bigg)
 \frac{\log x}{\sqrt q}
 +O(1)
\]
for general multiplicative functions, and $4/\sqrt{6}$ replaced by 1 for completely multiplicative functions. We finally let $x\to\infty$ first
with $\delta$ fixed, and then let $\delta\to0$. This proves the
proposition.
\end{proof}

\section{Completion of the proof}\label{sec:completion}

\begin{proof}[Proof of Theorem~\ref{thm:main}]
From \eqref{eq:decomposition},
Proposition~\ref{prop:sparse}, and Proposition~\ref{prop:two-large}, we get
\[
 |L_f(x,a/q)|
 \leq
 \bigg(\frac4{\sqrt6}+o(1)\bigg)
 \frac{\log x}{\sqrt q}
 +O(\log\log x).
\]
Since $\log\log x\ll\log q$ in this range, this proves
\eqref{eq:main} at $a/q$, and \eqref{eq:phase-error} returns us to
$\alpha$.
If $f$ is completely multiplicative, the second part of
Proposition~\ref{prop:two-large} gives us the same result with $4/\sqrt{6}$ replaced by 1.

The constant $4/\sqrt6$ is best possible by
\cite[Example~1.7]{GranvilleLamzouri2026}. To see that the constant $1$
is also best possible, we fix a primitive character $\chi\pmod{q_0}$ and
take $f=\overline{\chi}$ and $\alpha=a/q_0$, where $(a,q_0)=1$. Then
\[
 L_f(x,a/q_0)
 =
 \frac{\chi(a)\tau(\overline{\chi})}{q_0}\log x+O(1) \implies
 |L_f(x,a/q_0)|
 =
 (1+o(1))\frac{\log x}{\sqrt{q_0}}.
\]
\end{proof}

\end{document}